\documentclass{amsart}

\usepackage{amscd}
\usepackage{amssymb}
\theoremstyle{plain}

\newtheorem{thm}{Theorem}

\newtheorem{lem}[thm]{Lemma}

\theoremstyle{definition}
\newtheorem{df}[thm]{Definition}
\newtheorem{rem}[thm]{Remark}

\author{Michael Robinson} 

\title{Instability of an equilibrium with negative definite linearization}

\begin{document}
\begin{abstract}
A nonlinear parabolic differential equation is presented which has at
least one equilibrium.  This equilibrium is shown to have a negative
definite linearization, but a spectrum which includes zero.  An
elementary construction shows that the equilibrium is not stable.
\end{abstract}

% Keywords: stability, equilibrium, negative definite, spectrum
% MSC Classification: 37L15 (Primary); 35Q55 (Secondary)

\maketitle

\section{Introduction}
This note demonstrates that in infinite-dimensional settings, negative
definiteness of an equilibrium of a dynamical system is not sufficient
to ensure that the equilibrium is stable.  This is in stark contrast
to the situation in finite-dimensional settings, where negative
definiteness implies stability of the equilibrium.  (See
\cite{BoyceDiPrima}, for instance.)

The particular problem we study is the Cauchy problem
\begin{equation}
\label{pde}
\begin{cases}
\frac{\partial u(t,x)}{\partial t} = \Delta u(t,x) - 2 f(x) u(t,x) -
u^2(t,x)\\
u(0,x)=h(x) \in C^\infty(\mathbb{R}^n)\\
t>0,x\in \mathbb{R}^n \text{ for } n\ge 1,\\
\end{cases}
\end{equation}
where $f \in C_0^\infty(\mathbb{R}^n)$ is a positive function.  Since
the linear portion of the right side of \eqref{pde} is a sectorial
operator, we can use \eqref{pde} to define a nonlinear
semigroup. \cite{Henry} This turns \eqref{pde} into a dynamical
system, the behavior of which is largely controlled by its equilibria.
This problem evidently has as an equilibrium, $u(t,x) \equiv 0$ for
all $t,x$.  Depending on the exact choice of $f$, there may be other
equilibria, however, they will not concern us here.  The spectrum of
the equilibrium $u \equiv 0$ includes zero, even though the
linearization of \eqref{pde} about it is negative definite.  We show
this using an elementary construction akin to that of
\cite{Mazya_2005}.  Additionally, we show by a direct construction
that this equilibrium is not stable when $n=1$.

\section{Motivation}
The equation \eqref{pde} arises as a transformation of a related
equation, namely

\begin{equation}
\label{pde1}
\begin{cases}
\frac{\partial u(t,x)}{\partial t} = \Delta u(t,x) - u^2(t,x) + \phi(x)\\
u(0,x)=w(x) \in C^\infty(\mathbb{R}^n)\\
t>0,x\in \mathbb{R}^n \text{ for } n\ge 1,\\
\end{cases}
\end{equation}
with $\phi \in C_0^\infty(\mathbb{R}^n)$.  This equation describes a
reaction-diffusion equation \cite{FiedlerScheel}, or a diffusive
logistic population model with a spatially-varying carrying capacity.
The spatial inhomogeneity of $\phi$ makes the analysis of \eqref{pde1}
much more complicated than that of typical reaction-diffusion
equations.  The existence of the equilibria for \eqref{pde1} is a
fairly difficult problem, which depends delicately on $\phi$.  We will
not treat the existence of equilibria for \eqref{pde1} here, but
assume that $f$ is a positive equilibrium for \eqref{pde1}.  Then we
can look at the behavior of perturbations near $f$, for instance
\begin{eqnarray*}
\frac{\partial (f+u)}{\partial t} &=& \Delta (f+u) - (f+u)^2 + \phi\\
\frac{\partial u}{\partial t} &=& \Delta f+ \Delta u - f^2 -2fu -u^2 +
\phi\\
\frac{\partial u}{\partial t} &=& \Delta u -2fu -u^2,
\end{eqnarray*}
which is \eqref{pde}.  Notice that this transforms the equilibrium $f$
of \eqref{pde1} to the zero function in \eqref{pde}.  The situation of
\eqref{pde} is considerably easier to examine.

\section{Properties of the spectrum}
\label{stability_sec}
We need to linearize \eqref{pde} in order to examine the spectrum of
the equilibrium.  In doing so, we roughly follow the outline given
in \cite{Henry}.  Recall the following definition of the derivative
map in a Banach space:

\begin{df}
Suppose $R:B_1 \to B_2$ is a map from one Banach space to another.
The {\it derivative map} of $R$ at $u \in B_1$ is the unique linear
map $D:B_1 \to B_2$ such that for each sequence $\{h_n\}_{n=1}^\infty$
with $\|h_n\| \to 0$, 
\begin{equation*}
\lim_{n \to \infty} \left\| \frac{D(h_n)-R(u+h_n)+R(u)}{\|h_n\|}
\right\| = 0. 
\end{equation*}
Of course, such a map may not exist.  If it does, we say $R$ is {\it
  differentiable} at $u$.  The {\it linearization} $L$ of $R$ is the
  affine map given by the formula $L(h) = R(u) + D(h)$.
\end{df}

For this section, we shall work in the Hilbert space
$L^2(\mathbb{R}^n)$ with the usual norm (using the fact that $\Delta$
is densely defined wherever necessary).  The linearization of
\eqref{pde} at $u\equiv 0$ is easily computed to be
\begin{equation}
\label{linear_pde1}
\frac{\partial h(t,x)}{\partial t} = \Delta h(t,x) - 2 f(x) h(t,x).
\end{equation}

Suppose $h(x,t)=X(x)T(t)$, then we can separate variables in
\eqref{linear_pde1}, obtaining
\begin{eqnarray*}
T'(t)-\lambda T(t) = 0 \\
\Delta X(x) - (\lambda + 2 f(x))X(x) = 0.
\end{eqnarray*}
The separation constant $\lambda$ can be determined by examining the
eigenvalue problem
\begin{equation}
\label{schrod_eig}
(\Delta - 2 f(x)) X(x) = \lambda X(x),
\end{equation}
which is essentially the computation of the energy levels of a
Schr\"{o}dinger equation.  The operator $(\Delta - 2f)$ is a
Schr\"{o}dinger operator with potential $-2f$.  Due to its importance
in quantum mechanics, much is known about Schr\"{o}dinger
operators (see \cite{Mazya_2007} for a summary).

If $\Re(\lambda)<0$ over all of the eigenvalues $\lambda$ in
\eqref{schrod_eig}, we would normally conclude that $h \to 0$ as $t
\to \infty$, that $u \equiv 0$ is a stable equilibrium.  However, as
we shall see in Section \ref{instability_sec}, this is false.  The
cause of the instability is that although $\Re(\lambda)<0$ for all
eigenvalues, $\lambda=0$ is in the spectrum of the operator $(\Delta -
2f)$.

\begin{lem}
\label{selfadjoint_negdef_lem}
The spectrum of a self-adjoint, negative definite operator $T$ has
spectrum which is confined to the closed left half-plane $\{\lambda \in
\mathbb{C} | \Re (\lambda) \le 0 \}$.
\begin{proof} 
This is a standard argument (for instance, see \cite{Kreyszig}), which
we sketch briefly.  First, suppose $\lambda$ is an eigenvalue of $T$
with an eigenfunction $\psi$.  Then

\begin{equation*}
\lambda = \frac{\left<\psi,T \psi\right>}{\left<\psi,\psi\right>} = \frac{\left<T \psi,
  \psi\right>}{\left<\psi,\psi\right>} =\bar{\lambda} \le 0.
\end{equation*}
On the other hand, the Fredholm alternative (see \cite{Keener})
implies that $T-\lambda$ is surjective for $\lambda > 0$.

Finally, we note that for $\Re(\lambda) > 0$, $(T-\lambda)^{-1}$ is bounded:
\begin{eqnarray*}
\left<(T-\lambda)\psi,(T-\lambda)\psi\right>
&=&\left<T\psi,T\psi\right>-2\Re(\lambda)\left<\psi,T\psi\right>+|\lambda|^2\left<\psi,\psi\right>\\
&\ge&|\lambda|^2\left<\psi,\psi\right>,\\
\end{eqnarray*}
by the negative definiteness of $T$.  Hence, for $\Re(\lambda)>0$,
$(T-\lambda)$ has a bounded inverse.
\end{proof}
\end{lem}

\begin{lem}
\label{spec_sign_lem}
The self-adjoint operator $(\Delta - 2 f(x))$ is negative definite if and only if
$f > 0$ almost everywhere.  (See \cite{Mazya_2005} for a
generalization.)
\begin{proof}
It is well-known and easily shown that $(\Delta - 2f)$ is
self-adjoint.  See \cite{Kato_1951}, for example.  The
self-adjointness of $(\Delta - 2f)$ follows immediately from that of
$\Delta$.  It is also well-known that $\Delta$ is negative definite:
with zero boundary conditions, the divergence theorem gives
\begin{eqnarray*}
\left<u,\Delta u\right> &=& \int \bar{u} \Delta u dx \\
&=&-\int \nabla \bar{u} \cdot \nabla u dx < 0.\\
\end{eqnarray*}
So the only thing that will spoil the negative definiteness is $f$.
Suppose $f> 0$ almost everywhere, and $u \in L^2$.  Then
\begin{equation*}
\left<u,-2fu\right>=-2\int \bar{u} fu dx = -2\int f |u|^2 dx < 0.
\end{equation*}
On the other hand, suppose $A=\{x \in \mathbb{R}^n|f(x)\le 0\}$ has
positive measure.  Then let $u = 1_A$ and compute
\begin{equation*}
\left<u,(\Delta-2f)u\right>=\left<u,-2fu\right> = -2 \int \bar{u} f u
dx = -2 \int f |u|^2 dx \ge 0.
\end{equation*}
So we have that $(\Delta - 2 f)$ is not negative definite in that case.
\end{proof}
\end{lem}

\begin{lem}
\label{nonsingularity}
Suppose $f$ is a positive continuous function on $\mathbb{R}^n$.  Then
$(\Delta - 2f)$ is injective on $C_0^2(\mathbb{R}^n)$.
\begin{proof}
Let $u \in C_0^2(\mathbb{R}^n)$ satisfy $(\Delta - 2f)u=0$.  Let
$y=\sup_{x\in\mathbb{R}^n} u(x)$.  We claim that $y=0$.  Suppose the
contrary, that $y>0$.  Since $u\in C_0^2(\mathbb{R}^n)$, there is an
$R>0$ such that for all $\|x\|>R$, $u(x)<y$.  Thus $M=u^{-1}(\{y\})$
is compact.  By the maximum principle, there exists an $\epsilon>0$
such that the $\epsilon$-neighborhood of $M$, 
\begin{equation*}
M_\epsilon=\{x\in\mathbb{R}^n | \inf_{z\in M} \|z-x\| < \epsilon\}
\end{equation*}
has $\Delta u|(M_\epsilon-M)<0$.  On the other hand, $N=M_\epsilon
\cap u^{-1}((0,y))$ is an open set on which $u|N>0$ and $\Delta u |
N<0$.  But since $f$ is positive and $\Delta u = 2fu$, this is a
contradiction.  Similar reasoning leads to $\inf_{x\in\mathbb{R}^n}
u(x) = 0$, so in fact $u \equiv 0$.
\end{proof}
\end{lem}

Since $C_0^2(\mathbb{R}^n)$ is dense in $L^2(\mathbb{R}^n)$, this
implies that $\lambda=0$ is not an eigenvalue of $(\Delta-2f)$ over
$L^2(\mathbb{R}^n)$.

\begin{lem}
\label{supremum}
The spectrum of $(\Delta - 2f)$ includes zero when $f \in
C_0^\infty(\mathbb{R}^n)$ is a positive function.  (See
\cite{Mazya_2005} for the most general result of this kind.)
\begin{proof}
By Lemma \ref{nonsingularity} and the Fredholm alternative,
$(\Delta-2f)^{-1}$ exists.  We show that $(\Delta - 2f)^{-1}$ is not
bounded, by constructing a sequence $\{\psi_m\}$ such that
\begin{equation*}
\lim_{m \to \infty} \frac{\left<(\Delta-2f)\psi_m,(\Delta -
  2f)\psi_m\right>}{\left<\psi_m,\psi_m\right>} = 0.
\end{equation*}

Let $\psi_m$ be the function
\begin{equation*}
\psi_m(x)=\left(\frac{1}{2A_m\sqrt{\pi}}\right)^{n/2} e^{-\frac{\|x-B_m\|^2}{2A_m^2}},
\end{equation*}
where $A_m \in \mathbb{R}$ and $B_m \in \mathbb{R}^n$ are constructed
as follows.  Choose $A_m$ so that 
\begin{equation*}
\left<\Delta \psi_m, \Delta \psi_m \right> < \frac{1}{2m}
\end{equation*}
(that this is possible follows from an easy
computation).  
% Use MAPLE to do it.  Works out that <\Delta \psi_m,\Delta \psi_m> is roughly 1/A^4
Then select $B_m$ so that 
\begin{equation*}
\left<f \psi_m, f \psi_m\right> < \frac{1}{2m},
\end{equation*}
which is possible since $f \in C_0(\mathbb{R}^n)$.
Notice that $\left<\Delta \psi_m, \Delta \psi_m\right>$ is independent of $B_m$, so the
second choice does not interfere with the first.  Evidently
\begin{equation*}
\lim_{m \to \infty}\left<(\Delta-2f)\psi_m,(\Delta-2f)\psi_m\right> = 0,
\end{equation*}  
by the Schwarz inequality.  On the other hand, $\|\psi_m\|_2=1$ for
all $m$.  As a result, this shows that $(\Delta-2f)^{-1}$ is not
bounded.
\end{proof}
\end{lem}

As a result of Lemmas \ref{nonsingularity} and \ref{supremum}, we have
three things: (1) that that the spectrum is contained in the closed
left half plane, (2) the spectrum includes zero, and (3) zero is not
an eigenvalue. 

\section{Instability of the equilibrium}
\label{instability_sec}

Now we construct, for each $\epsilon>0$ and $1\le p < \infty$, an
$h_\epsilon \in C_C^\infty \cap L^p(\mathbb{R})$ such that
$\|h_\epsilon\|_p < \epsilon$ which if $u$ solves \eqref{pde} with
$h_\epsilon$ as its initial condition, then $\|u(t,\cdot)\|_p \to
\infty$.  In particular, this implies that $u\equiv 0$ is not a stable
equilibrium of \eqref{pde}.  We follow the general idea of the first
part of \cite{Fujita}.  (Additionally, \cite{Evans} contains a more
elementary discussion with a similar construction.)

\begin{df}
Let $H(t,x)=\frac{1}{\sqrt{4\pi t}} \exp\left(-
\frac{|x|^2}{4t}\right)$, which is the heat kernel.  Let
$v_\epsilon(s,x)=H(t-s+\epsilon,x)$ for fixed $t$ and $s<t$.  
\end{df}

\begin{rem}
Since $H$ is the heat kernel, $v_\epsilon$ will satisfy
$\frac{\partial v_\epsilon(s,x)}{\partial s} = - \Delta
v_\epsilon(s,x)$.  
\end{rem}

\begin{lem}
\label{fence_lem}
Suppose $u(t,x)\le 0$ satisfies \eqref{pde}, and $u(t,\cdot)\in
L^p(\mathbb{R})$ for each $t$.  Define 
\begin{equation}
\label{jepsilon}
J_\epsilon(s)=\int{v_\epsilon(s,x)u(s,x) dx}.
\end{equation}
Then $\frac{d J_\epsilon(s)}{ds} \le -(J_\epsilon(s))^2 -2
\|f\|_\infty J_\epsilon(s)$.
\begin{proof}
First of all, we observe that since $u\in L^p$,
  $v_\epsilon(s,\cdot)u(s,\cdot)$ is in $L^1(\mathbb{R})$ for each
$s<t$.

Now suppose we have a sequence $\{m_n\}$ of compactly supported
smooth functions with the following properties: \cite{LeeSmooth}
\begin{itemize}
\item $m_n \in C^\infty(\mathbb{R})$,
\item $m_n(x) \ge 0$ for all $x$, 
\item $\text{supp}(m_n)$ is contained in the interval $(-n-1,n+1)$, and
\item $m_n(x)=1$ for $|x| \le n$.
\end{itemize}
Then it follows that 
\begin{equation*}
J_\epsilon(s)=\lim_{n\rightarrow \infty} \int{v_\epsilon(s,x)u(s,x)m_n(x) dx}.
\end{equation*}

Now 
\begin{eqnarray*}
\frac{d}{ds} J_\epsilon(s) &=& \frac{d}{ds} \lim_{n\rightarrow \infty}
\int{v_\epsilon(s,x)u(s,x)m_n(x) dx} \\
&=& \lim_{h\rightarrow 0} \lim_{n\rightarrow \infty} \frac{1}{h}\int
(v_\epsilon(s+h,x)u(s+h,x)-v_\epsilon(s,x)u(s,x))m_n(x) dx.\\
\end{eqnarray*}
We'd like to exchange limits using uniform convergence.  To do this we
show that 
\begin{equation}
\label{big_lim}
\lim_{n\rightarrow \infty} \lim_{h\rightarrow 0} \frac{1}{h}\int
(v_\epsilon(s+h,x)u(s+h,x)-v_\epsilon(s,x)u(s,x))m_n(x) dx
\end{equation}
exists and the inner limit is uniform.  We show both together by a
little computation, using uniform convergence and LDCT:

\begin{eqnarray*}
&&\lim_{n\rightarrow \infty} \lim_{h\rightarrow 0} \frac{1}{h}\int
(v_\epsilon(s+h,x)u(s+h,x)-v_\epsilon(s,x)u(s,x))m_n(x) dx\\
&=&\lim_{n\rightarrow \infty} \int
\left(\frac{d}{ds}v_\epsilon(s,x)u(s,x)+v_\epsilon(s,x)\frac{d}{ds}u(s,x)\right)m_n(x) dx \\
&=&\lim_{n\rightarrow \infty} \int
(-\Delta v_\epsilon(s,x)u(s,x)+v_\epsilon(s,x)(\Delta u(s,x) -
u^2(s,x)- 2 f(x) u(x))m_n(x) dx \\
&=&\lim_{n\rightarrow \infty} \int
(-v_\epsilon(s,x)u^2(s,x)-2v_\epsilon(s,x)f(x)u(s,x))m_n(x) dx.
\end{eqnarray*}
Minkowski's inequality has that
\begin{equation*}
\int v_\epsilon u m_n dx \le \left(\int v_\epsilon m_n dx\right)^{1/2}
\left(\int v_\epsilon u^2 m_n dx \right)^{1/2},
\end{equation*}
since $v_\epsilon, m_n\ge 0$.  This gives that

\begin{eqnarray*}
&&\int
(-v_\epsilon(s,x)u^2(s,x)-2 v_\epsilon(s,x)f(x)u(s,x))m_n(x) dx\\
&\le& - \frac{(\int v_e u m_n dx)^2 }{\int v_\epsilon m_n dx} -
  2 \|f\|_\infty \int v_\epsilon u m_n dx\\
&\le& - \frac{\left(\int v_e u dx \right)^2}{\int v_\epsilon m_1 dx} -
  2 \|f\|_\infty J_\epsilon(s) < \infty,\\
\end{eqnarray*}
hence the inner limit of \eqref{big_lim} is uniform.  On the other
hand, 

\begin{eqnarray*}
&&\lim_{n\rightarrow \infty} \int
(-v_\epsilon(s,x)u^2(s,x)-2 v_\epsilon(s,x)f(x)u(s,x))m_n(x) dx\\
&\le&\lim_{n\rightarrow \infty} \left( - \frac{(\int v_e u m_n dx)^2
  }{\int v_\epsilon m_n dx} -
  2\|f\|_\infty \int v_\epsilon u m_n dx\right)\\
&\le&-(J_\epsilon(s))^2 - 2 \|f\|_\infty J_\epsilon(s)< \infty,\\
\end{eqnarray*}
so the double limit of \eqref{big_lim} exists.  Hence we conclude that 
the lemma is true.
\end{proof}
\end{lem}

\begin{lem}
\label{blowup_lem}
Suppose that for some $t_0>0$, 
\begin{equation*}
\int H(t_0,x) u(0,x) dx < - 2 \|f\|_\infty.
\end{equation*}
Then $\|u(t,\cdot)\|_p \to \infty$ for $1 \le p \le \infty$.
\begin{proof}
Note that 
\begin{eqnarray*}
J_\epsilon(0)&=&\int v_\epsilon(0,x) u(0,x) dx \\
&=& \int H(t+\epsilon,x) u(0,x) dx \\
& < & - 2\|f\|_\infty,
\end{eqnarray*}
since we may choose $\epsilon>0$ and $t$ such that $t+\epsilon=t_0$.
Thus Lemma \ref{fence_lem} implies that $J_\epsilon(s) \to -\infty$ by
elementary ODE theory.  \cite{BoyceDiPrima}

On the other hand,
\begin{eqnarray*}
|J_\epsilon(s)| \le \int | v_\epsilon(s,x) | |u(s,x)| dx &\le&
 \frac{1}{\sqrt{4\pi \epsilon}} \|u(s,\cdot)\|_1\\
&\le& \|u(s,\cdot)\|_\infty.\\
\end{eqnarray*}
So we have that $\|u(s,\cdot)\|_1$ and $\|u(s,\cdot)\|_\infty$ both
blow up.  Finally, 
\begin{eqnarray*}
\int | v_\epsilon(s,x) | |u(s,x)| dx &\le& \int |v_\epsilon| |u|^p
|u|^{1-p} dx \\
& \le & \frac{1}{\|u\|^{p-1}_\infty \sqrt{4\pi\epsilon}} \|u\|^p_p\\
&\le& \frac{1}{\sqrt{4\pi\epsilon}} \|u\|^p_p\\
\end{eqnarray*}
since $\|u(s,\cdot)\|_\infty \to \infty$.  Hence $\|u(s,\cdot)\|_p \to
\infty$.
\end{proof}
\end{lem}

Finally, we show that $u \equiv 0$ is unstable.  Let $\epsilon > 0$ be
given and $1 \le p < \infty$.  Take $h_\epsilon(0) \le -4 \|f\|_\infty$ to be
arbitrary.  We can construct $h_\epsilon \in L^p \cap L^\infty \cap
C^\infty_C(\mathbb{R})$ such that additionally $\|h_\epsilon\|_p < \epsilon$,
using the smooth Urysohn lemma.  \cite{LeeSmooth}  Then for
sufficiently small $t>0$,
\begin{equation*}
\int H(t,x) h_\epsilon(x) dx < - 2 \|f\|_\infty
\end{equation*}
by the fact that $\{H(1/n,\cdot)\}$ is a $\delta$-sequence as $n \to
\infty$.  Hence by Lemma \ref{blowup_lem}, if $u$ solves \eqref{pde}
with $h_\epsilon$ as initial condition, then $\|u\|_p \to \infty$.
(Note that this construction fails for $p=\infty$, since we cannot
ensure that both $h_\epsilon(0) \le - 4 \|f\|_\infty$ and
$\|h_\epsilon\|_\infty < \epsilon$.)

\section{Conclusions}
As a result of the previous two sections, we conclude that the
equilibrium $u \equiv 0$ has a real, negative eigenvalues (the set of
which may be negative), yet it is not stable.  That there exist
solutions which start near the equilibrium but blow up to $\infty$ in
any $p$-norm indicates that the equilibrium is actually rather
unstable.  On the other hand, this is precisely the kind of behavior
that is expected from a linearization whose spectrum contains zero.
This suggests that in infinite-dimensional settings one should be sure
to employ the entire spectrum to determine stability.

\bibliography{instability_bib}
\bibliographystyle{hplain}

\end{document}